\documentclass[a4paper]{article}

\usepackage{amsfonts}
\usepackage{amssymb}
\usepackage{latexsym}

\newtheorem{teo}{Theorem}[section]

\newtheorem{lem}[teo]{Lemma}

\newtheorem{prob}[teo]{Problem}
\newtheorem{conj}[teo]{Conjecture}

\newtheorem{defin}[teo]{Definition}
\newcommand{\dok}[1]{{\bf Proof.} $\:$ {#1} $\qquad \square$ \medskip}

\begin{document}


\title{$PI$-groups and $PI$-representations of groups}

\maketitle    

\begin{center}

\author{E. Aladova$^{1}$
      \and B. Plotkin$^{2}$}

 \smallskip
        {\small
                $^{1}$
          Department of Mathematics,
          Bar-Ilan University,
          52900, Ramat Gan, Israel

                {\it E-mail address:} aladovael@mail.ru
        }

\smallskip

 \smallskip
        {\small
               $^{2}$ Department of Mathematics,
          Hebrew University of Jerusalem,
          91904, Jerusalem, Israel

                {\it E-mail address:} plotkin@macs.biu.ac.il
        }
\end{center}

\begin{abstract}

It is well known that many famous Burnside-type problems have  positive
solutions for $PI$-groups and $PI$-algebras. In the present
article we also consider various Burnside-type problems for
$PI$-groups and $PI$-representations of groups.

\end{abstract}

\tableofcontents


\section{Introduction}\label{S_Int}

We start with a short historical background. The General Burnside
problem asks: {\it  Is a torsion group locally finite?}
 In 1964 E.S.~Golod obtained
negative solution of this problem: he constructed infinite
finitely generated residual finite torsion groups. His result
follows from the theorem of
Golod-Shafarevich~\cite{Golod-Shafarevich}. However,
C.~Procesi~\cite{Procesi} and A.~Tokarenko~\cite{Tokarenko}
obtained a positive solution of the Burnside Problem for
$PI$-groups: every periodic $PI$-group is locally finite.

The General Burnside problem gives rise to numerous questions of
Burnside-type which have positive answers in the case of
$PI$-groups and $PI$-algebras. Remind basic definitions.

Let $K$ be a field. Let $A$ be an associative $PI$-algebra with
unit over $K$, that is an algebra satisfying some non-trivial
polynomial identity.

\begin{defin}
A group $G$ is called a $PI$-group if there is a injective
homomorphism $\rho :G\to A^{*}$ of the group $G$ to the group
$A^{*}$ of invertible elements of the $PI$-algebra $A$.
\end{defin}

\begin{defin}
 An element $g$ of a group $G$ is called a nil-element if for every $x\in
G$ there is $n=n(g,x)$ such that $[[x,\underbrace{g],\dots
,g}_{n}]=1$.
\end{defin}

\begin{defin}
A group $G$ is called a nil-group if every its element is a
nil-element.
\end{defin}

\sloppy
\begin{defin}
 A group \ $G$ is called  Engel if it satisfies the identity \ \
 $[[x,\underbrace{y],\dots ,y}_{n}]=1$ for some $n$.
\end{defin}

Every locally nilpotent group is a nil-group, but the opposite statement in
general is not true: the theorem of
Golod-Shafarevich~\cite{Golod-Shafarevich} gives a negative
counterexample.
 For the case of  $PI$-groups
there is  a positive solution:
every nil-$PI$-group (and then Engel group) is locally nilpotent (see \cite{Platonov}, \cite{Plotkin_Notes}).

Moreover for Engel groups there is a long-standing conjecture:

\begin{conj}
An Engel group is not necessary locally nilpotent.
\end{conj}

This problem is still open: it is not known whether there exists a  non-locally
nilpotent Engel group.

For algebraic algebras the following Burnside-type problem posted by A.Kurosh is well-known.
 Remind the necessary definitions.

\begin{defin}
An algebra $A$ is called  algebraic  if every element
$a\in A$ satisfies some polynomial identity $f(x)=0$.
\end{defin}

\begin{defin}
An algebra $A$ is called  locally finite if every finite set of
elements of the algebra $A$ generates a finite dimensional
subalgebra.
\end{defin}

Every locally finite algebra is an algebraic algebra. The problem
of Kurosh asks: {\it Is every algebraic algebra locally finite?}
In general it is not true, but I.~Kaplansky and A.~Shirshov ~\cite{Jacobson} give
positive answer in the case of $PI$-algebras:  every algebraic
$PI$-algebra is locally finite.

We consider a variant of Kurosh's problem for groups. Below are the
necessary definitions.

\begin{defin}
An element $g$ of a group $G$ is called an algebraic element if
for every $x\in G$ the subgroup generated by the all elements of
the form $[[x,\underbrace{g],\dots ,g}_{n}],$ $n\in \mathbb
N$ is finitely generated.
\end{defin}

\begin{defin}
A group $G$ is called algebraic if every its element is an
algebraic element.
\end{defin}


\begin{defin}
A group $G$ is called locally Noetherian if every its finitely
generated subgroup is Noetherian.
\end{defin}

It is obvious that every locally Noetherian group is algebraic. We
consider the question when the inverse statement is true.
In Section 2 we show that there is the positive solution of this problem for $PI$-groups.

\begin{teo}\label{Th_AlgPI-gr}
Every algebraic $PI$-group is locally Noetherian.
\end{teo}

In fact we will prove that every finitely generated algebraic
$PI$-group is a Hirsch (polycyclic-by-finite)  group.

The  other our main result deals with  $PI$-representations of
groups. 

 Let $K$ be a field, $V$ be a $K$-module. Let $G$ be a group and
let $(V,G)$ be a representation of the group $G$ in the $K$-module
$V$. So we have a homomorphism $\rho$ of the group $G$ to the
group of automorphisms of the $K$-module $V$:
$$
\rho : G\to \mbox{Aut}V, \ \ g\to g^{\rho}.
$$

We can say that the group $G$ acts on the module $V$ by the rule:
$$
(v,g)\to v\circ g=g^{\rho}(v),
$$
 for all $v\in V$, $g\in G$.
 Note
that the group algebra $KG$ also acts on the module $V$.

Let $F$ be a free group of countable rank with free generators
$x_1,x_2,\dots $, and let $KF$ be the group algebra of $F$. Let
$u(x_1,\dots ,x_n)$ be an element of $KF$.

\begin{defin}
We say that a representation $(V,G)$ satisfies an identity $y\circ
u(x_1,\dots ,x_n)\equiv 0$  if for all $v\in V$ and all $g_i\in G$
we have $v\circ u(g_1,\dots ,g_n)=0$. 
\end{defin}

\begin{defin}
An element $g\in G$ is called a unipotent element of a
representation $(V,G)$ if  there is $n=n(g)$ such that $x\circ
(g-1)^{n}\equiv 0$ for every $x\in V$.
\end{defin}

\begin{defin}
A representation $(V,G)$ is called a unipotent representation if
every $g\in G$ is a unipotent element of $(V,G)$.
\end{defin}

\begin{defin}
A representation $(V,G)$ is called a unitriangular if it satisfies
the identity $x\circ (y_1-1)\dots (y_n-1)\equiv 0$.
\end{defin}

\begin{defin}
A representation $(V,G)$ is called a locally unitriangular if for
every finitely generated subgroup $H$ of the group $G$ the
subrepresentation $(V,H)$ is unitriangular.
\end{defin}

Every locally unitriangular representation is unipotent.
There is the following Burnside-type problem for unipotent
representations:

\begin{prob}\label{Pr_UnipRepr}
Is every unipotent representation locally unitriangular?
\end{prob}

In general, it is not true.
We prove (Section 3) the positive solution of this problem for
$PI$-representations of groups.

Let $(V,G)$ be a representation of the group $G$ in the $K$-module
$V$. Let $\overline G=G/\mbox{Ker}(V,G)$ and $(V,\overline G)$ be
the faithful representation corresponding to the representation
$(V,G)$.

\begin{defin}\label{Def_PI-repr}
A representation $(V,G)$ is called a $PI$-representation if the
linear span $\langle \overline G\rangle $ of the group $\overline
G$ in the algebra $\mbox{End }V$ is a PI-algebra.
\end{defin}

We have the following Theorem (see Section 3)

\begin{teo}\label{Th_UnipPI-Repr}
Every unipotent $PI$-representation is locally unitriangular.
\end{teo}

In the frameworks  of  studying the unipotent representations we also
consider the following problem:

\begin{prob}\label{}
Is there a unipotent radical
 for a representation of a group?
\end{prob}

It is well-known that every group has a locally nilpotent radical
(Hirsch-Plotkin radical) that is unique maximal normal locally
nilpotent subgroup~\cite{Hirsch}, \cite{Plotkin_LocNilpRad}, but
it is not true for locally solvable radical, since there exist groups
without locally solvable radical (see results of G.~Baumslag,
L.~Kovach, B.~Neumann, V.~Mikaelian). However S.~Pikhtilkov proved
that every $PI$-group has a locally solvable
radical~\cite{Pihtilkov}.

Let $(V,G)$ be a representation of the group $G$ in the module
$V$.

\begin{defin}\label{Def_UnipRad}
The unique maximal normal subgroup $H$ of the group $G$, such that
the subrepresentation $(V,H)$ is locally unitriangular, is called
a locally unitriangular radical of the representation $(V,G)$.
\end{defin}

\begin{defin}\label{Def_UnipRad}
The unique maximal normal subgroup $H$ of the group $G$, such that
the subrepresentation $(V,H)$ is unitriangular, is called a
unipotent radical of the representation $(V,G)$.
\end{defin}

It is known~\cite{Plotkin_Autom} that for every representation of
a group there is a locally unitriangular radical. We have the
following

\begin{teo}\label{Th_UnipRad}
There exists a unipotent radical for a $PI$-representation and it 
coincides with the locally unitriangular radical.
\end{teo}

\noindent
which is an immediate corollary of Theorem \ref{Th_UnipPI-Repr}.


{\bf Acknowledgement.} The first author is supported by the Israel
Science Foundation  grant number 1178/06.

\section{Algebraic $PI$-groups}\label{S_Alg_PI-gr}

In this section we will prove the theorem~\ref{Th_AlgPI-gr} which
states that every algebraic $PI$-group is locally Noetherian.

\bigskip
{\bf Proof of the theorem~\ref{Th_AlgPI-gr}.} Let $G$ be an
algebraic finitely generated $PI$-group. Since $G$ is a $PI$-group
then there exists a $PI$-algebra $A$ such that the group $G$ is a
subgroup of the group of invertible elements of the algebra $A$.
Let $R(A)$ be the Levitzky radical of the algebra $A$.
%
%
Consider the algebra $A/R(A)$. The group $G$ acts on $A/R(G)$ and
the kernel of this action is $H=G\cap (1+R(A))$, moreover the
group $H$ is a locally nilpotent subgroup of the group
$G$~\cite{Plotkin_Notes}.

Now consider the group $\overline G=G/H$. It is
known~\cite{Jacobson} that there is an embedding $A/R(A)\to
M_n(K)$, where $M_n(K)$ is the matrix algebra of dimension $n$ and
$K$ is a commutative ring with unit which is a Cartesian sum of
fields. So the group $\overline G$ is a subgroup of the group
$GL_n(K)$ of invertible elements of the algebra $M_n(K)$.

Since the group $G$ is finitely generated algebraic group then the
group $\overline G$ is also finitely generated and algebraic. Let
$\overline G=\langle \overline g_1,\dots ,\overline g_m \rangle $,
$\overline g_i\in GL_n(K)$, $i=1,2,\dots ,m$. Let $\overline
g_i=(\alpha _{st}^{i})$ and let $S$ be a set of allelements
$\alpha _{st}^{i}$  such that  $\alpha _{st}^{i}\in K$, $i=1,2,\dots m$
and  $s,t=1,2,\dots n$. Note that the set $S$ is finite.

Let $K_0$ be a subring of the ring $K$ generated by the set $S$.
Since the ring $K_0$ is a finitely generated commutative ring then it is
Noetherian. From the theorem of Lasker~\cite{Zariski_Samuel}
follows that the ring $K_0$ has a finitely many prime ideals
$I_{\alpha }$ such that $\bigcap_\alpha I_{\alpha}=0$ and
$K_0/I_{\alpha }$ are fields. Then using the Remak's theorem we
have that $K_0$ is a Cartesian sum of the fields $K_0/I_{\alpha
}$.

Let $M_n(K_0)$ be the algebra of matrices over $K_0$. Then
$\overline G\subset M_{n}(K_0)$ and $M_n(K_0)$ is a subalgebra of
the algebra $\bigoplus_{\alpha} M_n(K_0/I_{\alpha})$, where
$\alpha\in \mathbb N$ and $\alpha< \infty$.

So the group $\overline G$ is embedded into a direct product of
the finite number of groups of matrices over fields:
$$
\overline G \to \prod_{\alpha}GL_n(K_0/I_{\alpha}), \ \alpha
<\infty.
$$

According to the well known Tits alternative \cite{Tits}  {\it
a finitely generated linear group either contains a non-abelian
free group or has a solvable subgroup of finite index.}

Since the group $\overline G$ is algebraic then every its subgroup
is algebraic, but a free group is not algebraic. So the group
$\overline G$ contains a solvable subgroup $\overline G_0$ of
finite index. It was proved  
 that every
locally soluble algebraic group is locally Noetherian and
consequently locally polycyclic (see ~\cite{Plotkin_Autom}).
%
 Thus the subgroup $\overline G_0$
is finitely generated  and locally polycyclic, 
 so it is a polycyclic group. 
 Moreover since  $\overline G_0$
is a subgroup of a finite index in $G$ then the group $\overline G$ is a
Hirsch group.

Remind that $\overline G=G/H$ and $H$ is a locally nilpotent
subgroup of the group $G$. If $\overline G_0=G_0/H$, then the
group $G_0$ is an extension of the locally nilpotent group $H$ by
the solvable group $\overline G_0$ and, hence, $G_0$ is algebraic. So the
group $G_0$ is locally Noetherian
and consequently it is a polycyclic group. Moreover $G_0$ is a
subgroup of a finite index in $G$ since $\overline G_0$ is a
subgroup of finite index in $\overline G$.
%
So in the group $G$ there exists the polycyclic subgroup $G_0$ of
 finite index and then $G$ is a Hirsch group. But every Hirsch group
is Noetherian. Thus $G$ is a Noetherian group.

Consequently every non-finitely generated algebraic $PI$-group is
locally Noetherian. Theorem is proved. $\qquad \square$

\section{Unipotent $PI$-representations of groups}\label{S5_Linearization}

%
%
In this section we will prove the theorem~\ref{Th_UnipPI-Repr}.
This theorem is a generalization of the well-known theorem of
Kolchin~\cite{Kolchin} which states:

\begin{teo}\label{Th_Kolchin}
Let $G$ be a linear group, $G\le GL_{n}(K)$, $K$ is a field. Let $(K^{n},G)$ be an
unipotent representation of the group $G$. Then the representation
$(K^{n},G)$ is unitriangular. $\qquad \square$
\end{teo}


Note that the representation $(K^{n},G)$ from this theorem
 is a $PI$-representation since the group $G$ is embedded into the
algebra $M_{n}(K)$
which satisfies the
standard identity of Amitsur-Levitzky.



To prove the theorem~\ref{Th_UnipPI-Repr} we need two lemmas.

\begin{lem}\label{Lem_5_1}
Let $(V,G)$ be a representation of a group $G$ and let $G$ is
generated by the set $M$. Let $\widehat{h}=(h_1,\dots , h_n)$ be a
sequence from the set $M^{n}$. If every $\widehat{h}$ satisfies
the equation $x\circ (h_1-1)\dots (h_n-1)=0$ for all $x\in V$ then
the representation $(V,G)$ satisfies the identity $x\circ
(y_1-1)\dots (y_n-1)\equiv 0$.
\end{lem}

\dok{
 The representation $(V,G)$ satisfies the identity
 $x\circ (y_1-1)\dots (y_n-1)\equiv 0$ if and only if the module
 $V$ has a $G$-invariant series of submodules of the length $n$
  $$
   0=V_n\subset V_{n-1}\subset \dots \subset V_1\subset V_0=V,
  $$
 such that the group $G$ acts trivially on the factors.

 We prove the lemma by induction on $n$.
 Let $n=1$ then for every $h\in M$ we have
 $$
  x\circ (h-1)=0 \mbox{ or } x\circ h=x,\  \mbox{for all } x\in V.
 $$

 For all $h_1,h_2$ from the set $M$ we have
 $$
  x\circ (h_1h_2)=(x\circ h_1)\circ h_2=x\circ h_2=x,
 $$
%
 then the group $G$ acts trivially on
 $V$ and we have the following series:
 $$
  0=V_1\subset V_0=V.
 $$

 Let now the statement of the lemma is hold for all positive
 integer less or equal $(n-1)$.
 Let $(h_1,\dots ,h_{n})$ be an arbitrary sequence in $M^{n}$ and
 let for every $x\in V$ we have
 $$
 x\circ (h_1-1)\dots (h_{n}-1)=0.
 $$
 Let $V_{n-1}$ be a linear span of all elements of the form $x\circ (h_1-1)\dots (h_{n-1}-1)$.
 The element $(h_{n}-1)$ annihilates the submodule
 $V_{n-1}$, since this element annihilates all generators of $V_{n-1}$
 and for every $v\in V_{n-1}$ we have $v\circ (h_n-1)=0$. Since $h_n$
 is an arbitrary element in $M$ then we have $v\circ (g-1)=0$ for all $g\in G$ and for all $v\in V_{n-1}$.
 So the group $G$ acts trivially on the module $V_{n-1}$.

 Consider the representation $(V/V_{n-1},G)$. For all
 $h_1,\dots ,h_{n-1}\in M$  and for all $x\in V$ we have
  $$
   x\circ (g_1-1)\dots (g_{n-1}-1)=0.
  $$
 Using the inductive assumption we have that the representation $(V/V_{n-1},G)$
 satisfies the identity
 $$
  x\circ (y_1-1)\dots (y_{n-1}-1)\equiv 0.
 $$
 It means that there is the following series of submodules of the
 module $V/V_{n-1}$:
  $$
   0=V_{n-1}/V_{n-1}\subset V_{n-2}/V_{n-1}\subset \dots \subset V_{1}/V_{n-1}\subset
   V/V_{n-1},
  $$
 such that the group $G$ acts trivially in the factors.

 Now consider the following series of submodules of the
 module $V$:
   $$
   0=V_{n}\subset V_{n-1}\subset \dots \subset V_{1}\subset V.
  $$
 The group $G$ acts trivially in the factors of this series. Thus the representation
 $(V,G)$ satisfies the identity
 $$
 x\circ (y_1-1)\dots (y_n-1)\equiv 0.
 $$
 Lemma is proved.
} 

Let $(V,G)$ be a $PI$-representation. Let $(V,\overline G)$ be a
faithful representation corresponding to the representation
$(V,G)$ and let $A=\langle \overline G \rangle $ be the linear
span of the group $\overline G$ in the algebra $\mbox{End }V$.
Note that $A$ is a $PI$-algebra.
 We can
consider the regular action of the group $\overline G$ in the
algebra $A$ assuming that the group $\overline G$ is embedded in
$A$ and so we have the representation $(A,\overline G)$ of the
group $\overline G$.

\begin{lem}\label{Lem_5_2}
Let $A_1$ be a nilpotent ideal of the algebra $A$ and let the
representation $(A/A_1,\overline G)$ satisfy the identity
$x\circ (y_1-1)\dots (y_{n}-1)\equiv 0$. Then the representations
$(A,\overline G)$, $(V,\overline G)$  and (V,G) are unitriangular.
\end{lem}

\dok{
 Let $\overline g_1,\dots ,\overline g_n$ be arbitrary elements of the group $\overline G$.
 Since the representation $(A/A_1,\overline G)$ satisfies the identity $x\circ
(y_1-1)\dots (y_{n}-1)\equiv 0$ then
 the element
 $(\overline g_1-1)\dots (\overline g_{n}-1)$ belongs to the ideal $A_1$.

 Let $A_1$ is
 nilpotent ideal of class nilpotency $m$ and let the elements
 $$
 \begin{array}{c}
 (\overline g_{11}-1)(\overline g_{12}-1)\dots  (\overline g_{1n}-1),\\
 (\overline g_{21}-1)(\overline g_{22}-1)\dots (\overline g_{2n}-1),\\
 \dots \\
 (\overline g_{m1}-1)(\overline g_{m2}-1)\dots (\overline g_{mn}-1)
 \end{array}
 $$
 belong to the ideal $A_1$, where $\overline g_{11},\dots ,\overline g_{mn}$
 are arbitrary elements of the group $\overline G$.
 The product of these elements is equal to zero:
 $$
 (\overline g_{11}-1)\dots (\overline g_{mn}-1)=0.
 $$
 Then the representations $(A,\overline G)$ and $(V,\overline G)$ satisfy the identity
 $$
  x\circ (y_{11}-1)\dots (y_{mn}-1)\equiv 0.
 $$
 So the representations $(A,\overline G)$ and $(V,\overline G)$ are unitriangular.


For varieties of representations of groups we have the following
invariant description~\cite{PV}:
{\it A class of representations of groups forms a variety of
representations of groups if and only it is closed under
subrepresentations, homomorphic images, Cartesian
products and 
saturations.}

Remind that a class of representations of groups $\mathfrak X$ is
closed under the saturation if the following condition is true: if
a representation $(V,H)$ lies in $\mathfrak X$ then all
representations $(V,G)$ such that $(V,H)$ is a right epimorphic
image of the representation $(V,G)$ also belong to $\mathfrak X$.

Thus the representation $(V,G)$ is also unitriangular.
} 

{\bf Remark.} For a $PI$-representation we have the following
properties: if a $PI$-representation $(V,G)$ satisfies the
identity
 $$
  x\circ (y_{1}-1)\dots (y_{n}-1)\equiv 0,
 $$
then the representation $(A,\overline G)$ also satisfies this
identity.

Indeed, let $(V,G)$ be a $PI$-representation, let $(V,\overline
G)$ be a faithul representation corresponding to $(V,G)$ and let
$A=\langle \overline G \rangle$. So we have also the faithful
representation $(V,A)$ of the algebra $A$ in the module $V$.

Let the $PI$-representation $(V,G)$ satisfies the identity
 $$
  x\circ (y_{1}-1)\dots (y_{n}-1)\equiv 0.
 $$
Then the faithful representation $(V,\overline G)$ also satisfies
this identity. So for every $x\in V$ the element $(\overline
g_1-1)\dots (\overline g_n-1)\in K\overline G$
 acts as zero. Since the representation $(V,A)$ is faithful we
 have $(\overline g_1-1)\dots (\overline g_n-1)=0$.
So for every $a\in A$ we have
 $$
 a\cdot (\overline g_1-1)\dots (\overline g_n-1)=0,
 $$
 and the regular representation $(A,\overline G)$ satisfies the identity
 $$
  x\circ (y_{1}-1)\dots (y_{n}-1)\equiv 0.
 $$

 In a similar manner, we can show that if a $PI$-representation $(V,G)$ satisfies the
 identity
 $$
  x\circ (g-1)^n=0,
 $$
 then 
 the representation $(A,\overline G)$ also satisfies this identity.
 Thus, if $g$ is a
 unipotent element of the representation $(V,G)$ then it is a
 unipotent element of the representation $(A,\overline G)$.

Now we can prove the theorem \ref{Th_UnipPI-Repr}.

\bigskip
{\bf Proof of the theorem.}
 Let $(V,G)$  be a unipotent $PI$-representation, let $(V,\overline G)$ be the faithful
 representation corresponding to the representation $(V,G)$ and let $A=\langle \overline G \rangle$ be
 a linear span of the group $\overline G$ in the algebra $\mbox{End }V$ . We can note that the regular
 representation $(A,\overline G)$ also is unipotent.

 Assume that the group $\overline G$ is finitely generated and let $\overline G=\mbox{gr}\langle M\rangle$,
 where $M$ is a finite set. We will prove that the representation
 $(V,\overline G)$ is unitriangular and hence the representation $(V,G)$ is also unitriangular~\cite{PV}.

 From a general structure theory \cite{Jacobson} it follows that the $PI$-algebra $A$ has
 a series of ideals
 $$
  0=A_0\subset A_1\subset A_2\subset A,
 $$
 where $A_1$ is a sum of nilpotent ideals of $A$,
 $A_2$ is the Levitzky radical of $A$,
%
%
 and the group $\overline G$ acts on the factors $A/A_1$ and $A/A_2$.

 We also have a series of subgroups of the group $\overline G$:
 $$
  \overline 1=\overline G_0\subset \overline G_1\subset \overline G_2\subset \overline G,
 $$
 where $\overline G_1$ is a kernel of the action $\overline G$ in $A/A_1$ and $\overline G_2$
 is a kernel of the action $\overline G$ in $A/A_2$.

Then,
 From a general theory \cite{Jacobson}
 it also follows that there is an
 injective homomorphism
 $$
  \tau :A/A_2\to M_{n}(R),
 $$
 where $M_{n}(R)$ is an algebra of matrices over a commutative ring
 $R$ with unit and the given ring $K$ is contained in $R$ (see \cite{Jacobson}). 
 We can consider the algebra
 $M_{n}(R)$ as an algebra over $K$, so the homomorphism $\tau$
 is a homomorphism of algebras over $K$. Moreover we have the
 representation $(R^{n},M_{n}(R))$.

The representation $(A/A_2,\overline G/\overline G_2)$ is unipotent
 as a homomorphic image of the unipotent representation $(A,\overline G)$,
 and an element $\overline g\in \overline G$ is a unipotent element of the representation $(A,\overline G)$
 if the element $(\overline g-\overline 1)$ is  a nilpotent element of the algebra $A$.

 A group $(\overline G/\overline G_2)^{\tau}$ is a subgroup of the group $GL_{n}(R)$ of all invertible
 elements of the algebra $M_{n}(R)$. For each element $\widetilde g\in
 G/G_2$ the element $(\widetilde g-\widetilde 1)$ from $A/A_2$ is nilpotent
 since the representation $(A/A_2,\overline G/\overline G_2)$ is unipotent. The
 image of this element in the algebra $M_{n}(R)$ is also
 nilpotent. Thus the representation $(R^{n}, (\overline G/\overline G_2)^{\tau})$ is
 a unipotent representation.
 According to the Kolchin's theorem $(R^{n},(\overline G/\overline G_2)^{\tau})$
 is a unitriangular representation.

 Well known theorem of Kaloujnine~\cite{Kaloujnine} states that if
 a representation $(V,\overline G)$ of a group $\overline G$ is faithful and $n$-unitriangular
 then the group $\overline G$ is a nilpotent group of the  nilpotency class $(n-1)$.

 Consider the representation $(R^{n}, \overline G/\overline G_2)$.  This
 representation is faithful since the regular representation
 $(A/A_2,\overline G/\overline G_2)$ and the representation  $(R^{n}, A/A_2)$ are faithful.
 Moreover the representation
 $(R^{n},(\overline G/\overline G_2)^{\tau})$ is unitriangular and therefore the representation
 $(R^{n},\overline G/\overline G_2)$ is also unitriangular.
 So the group $\overline G/\overline G_2$ is nilpotent.

 Since the
 representation $(A/A_2,\overline G/\overline G_2)$ is unipotent and the group
 $\overline G/\overline G_2$ is finitely generated nilpotent group
  then from \cite{Plotkin_Autom} 
 it follows that this representation is unitriangular.
 So the representation $(A/A_2,\overline G)$ is also unitriangular.

 Note that  $A_2/A_1$ is a nilpotent ideal of the algebra $A/A_1$ and $(A/A_1)/(A_2/A_1)\cong
 A/A_2$. Using  the lemma \ref{Lem_5_2}  we have that the representation
 $(A/A_1,\overline G)$ is unitriangular too.

 Let the representation $(A/A_1,\overline G)$ satisfy the identity
 $$
  x\circ (y_1-1)\dots (y_m-1)\equiv 0.
 $$
 Remind that we consider the finitely generated group $\overline G$ with a
 generating set $M$. Let $\widehat{g}=(\overline g_1,\dots ,\overline g_m)$ be a
 sequence of elements from the set $M^{m}$. Then the element
 $(\overline g_1-\overline 1)\dots (\overline g_{m}-\overline 1)$ lies in $A_1$.
 Since the set of different sequences $\widehat{g}$ is finite then there is a nilpotent
 ideal $A'_{1}\subset A_{1}$ such that for every $\widehat{g}$
 we have $(\overline g_1-\overline 1)\dots (\overline g_{m}-\overline 1)\in A'_{1}$ .
 Since $A'_{1}$ is an ideal of the algebra $A$ then it is closed
 under the regular action of the group $\overline G$ and we can consider the
 representation $(A/A'_{1},\overline G)$.
 According to the lemma \ref{Lem_5_1} the representation $(A/A'_{1},\overline
 G)$  satisfies the identity
 $$
  x\circ (y_1-1)\dots (y_{m}-1)\equiv 0.
 $$
 So this representation is unitriangular. Using the lemma
 \ref{Lem_5_2} we conclude that the representation $(V,G)$ is also
 unitriangular.


 For a non-finitely generated  group $G$ we have that
 the representation $(V,G)$ is locally unitriangular.
 Theorem is proved. $\qquad \square$

Theorem~\ref{Th_UnipPI-Repr} implies 
 Theorem~\ref{Th_UnipRad} which states that there exists a unipotent
radical for a $PI$-representation and it is coincides with the
locally unitriangular radical as a straightforward corollary.


\begin{thebibliography}{99}





\bibitem{Golod-Shafarevich}
E.S.~Golod, I.R.~Safarevich, On the class field tower, (Russian),
Izv. Akad. Nauk USSR Ser. Math., {\bf 28}, (1964), p.~261--272.



\bibitem{Hirsch}
K.A. Hirsch, \"Uber lokal-nilpotente Gruppen, Math.Z. {\bf 63},
(1955), 290--294.


\bibitem{Jacobson}
N.~Jacobson, Structure of rings, Providence: Amer. Math. Soc.,
(1964), 299~pp.

\bibitem{Kaloujnine}
L.~Kaloujnine, \" Uber gewisse Beziehungen zwischen einer Gruppe
und ihren Automorphismen. (German), Beriliner Math. Tagung,
(1953), pp. 164--172.



\bibitem{Kolchin}
E.R.~Kolchin, On certain concepts in the theory of algebraic
matrix groups, Ann. of Math., {\bf 49}, (1948), p. 774--789.






\bibitem{Novikov_Adian}
P.S.~Novikov, S.I.~Adian, On infinite periodic groups, I,II,III,
(Russian), Izv. Akad. Nauk USSR Ser. Math. {\bf 32}, (1968),
p.~212--244, p.~251--524, p.~709--731.


\bibitem{Pihtilkov}
S.A. Pikhtilkov, On the prime radical of PI-representable groups.
(Russian. Russian summary) Mat. Zametki {\bf 72} (2002), no. 5
739-744; English transl. im Math. Notes {\bf 72} (2002), no.5-6
682--686.


\bibitem{Platonov}
V.P.~Platonov, Engel elements and the radical in $PI$-algebras and
topological groups, (Russian) Doklady Akademii Nauk SSSR, {\bf
161}, N~2, (1965), p.~289--291. (English transl.: Sov. Math.
Dokl., {\bf 6}, (1965), p.~412--415.)


\bibitem{Plotkin_Autom}
B.I.~Plotkin, Group of automorphisms of algebraic system.
Translated from the Russian by K.A.~Hirsch. Wolters-Noordhoff
Publishing, Groningen, (1972), xviii+502~pp.


\bibitem{Plotkin_Notes}
B.I.~Plotkin, Notes on Engel groups and Engel elements in groups.
Some generalizations, (English, Russian summary), Izv. Ural. Gos.
Univ. Mat. Mekh., {\bf 36}, (2005), no.~7, p.~153--166,
p.~192--193.






\bibitem{Plotkin_LocNilpRad}
B.I. Plotkin, On  some criteria of locally nilpotet groups, Uspehi
Mat. Nauk, T. IX, {\bf 3(61)} (1954), 181-186. English transl. in
Amer. Math. Soc. Transl. (2) {\bf 17} (1961), 1--7.






\bibitem{PV}
B.I.~Plotkin, S.M.~Vovsi, Varieties of group representations.
General theory, connections and applications, (Rusian),
''Zinatne'', Riga, (1983), 339~pp.

\bibitem{Procesi}
C.~Procesi, The Burnside problem, J. Algebra, {\bf 4}, (1966), p.~421--425. 


\bibitem{Tits}
J.~Tits, Free subgroups in linear groups,J. Algebra, {\bf 20},
(1972), p.~250--270.

\bibitem{Tokarenko}
A.~Tokarenko, About linear groups over rings, Sibirsk. Mat. Z., {\bf 9},
(1968), no.~4, p.~951--959. 



\bibitem{Zariski_Samuel}

O. Zariski, P. Samuel, Commutetive algebra, Vol 1. With the
cooperation of I. S. Cohen. The University Series in Higher
Mathematics. D.Van Nostrand Company, Inc., Princeton, New Jersey,
(1958), xi+329 pp.














\end{thebibliography}
\end{document}